# COMPLETE ALGEBRAIC CLASSIFICATION OF THE HYPERSURFACES OF THE MAXIMUM INACCURACIES OF AN INDIRECTLY MEASURABLE VARIABLE

**Yordan Epitropov**[a], **Kiril Kolikov**[a], **Radka Koleva**[b]

[a]Plovdiv University "P. Hilendarski", Plovdiv, Bulgaria

[b]University of Food Technologies, Plovdiv, Bulgaria

Corresponding Author: Radka Koleva[b], r.p.koleva@gmail.com

[b]University of Food Technologies, 26 Maritsa Blvd, 4000 Plovdiv, Bulgaria

**Abstract**

Let an indirectly measurable variable $Y$ be represented as a function of a finite number of directly measurable variables $X_1, X_2, ..., X_n$. In our previous researches we: 1) represented the maximum inaccuracies of $Y$ in first degree of approximation as linear functions of the inaccuracies of $X_1, X_2, ..., X_n$; 2) defined the spaces of the maximum inaccuracies and we defined a dimensionless scale for quality (accuracy) evaluation of an experiment in them; 3) introduced the maximum inaccuracies in second degree of approximation.

In the current paper we prove that the maximum inaccuracies of $Y$ in second degree of approximation are quadrics of the inaccuracies of $X_1, X_2, ..., X_n$ and that these forms describe certain types of quadric hypersurfaces of parabolic class. Moreover: 1) we give a complete algebraic classification of these hypersurfaces; 2) we define a dimensionless scale for quality (accuracy) evaluation of the experiment given the maximum inaccuracies in second degree of approximation.

**1. Introduction.** Let $Y$ be an indirectly measurable variable which depends on $n$ directly measurable variables which are modelled (with the help of measuring instruments and methods) by $n$ real independent variables $X_1, X_2, ..., X_n$. Let us denote with $f$ the real function with arguments $X_i$ ($i = 1, 2, ..., n$) with the help of which one can represent $Y$, i.e. $Y = f(X_1, X_2, ..., X_n)$.

There are different methods [1, 2, 3, 4] for determining the inaccuracy (error) when determining the value of $Y$ in a given experiment. In [5, 6, 7] we researched the maximum absolute and maximum relative inaccuracy of the indirectly measurable variable $Y = f(X_1, X_2,..., X_n)$ in first degree of approximation. We suggested a new analytical principle for representation of the maximum inaccuracies of $Y$ as linear functions, respectively, of the absolute and relative inaccuracies of the directly measurable variables $X_i$ ($i = 1, 2,..., n$). We introduced spaces of the inaccuracies of $Y$. We defined sample planes of the ideal perfectly accurate experiment and using them a universal numerical characteristic – a dimensionless scale for quality (accuracy) evaluation of the experiment. With the help of the scale it can be shown how an experiment can be improved in order to reduce the measurement inaccuracy. Moreover, in [8] we introduced maximum absolute and relative inaccuracies of second order of $Y$.

In the current paper we show that the maximum inaccuracies of $Y$ in second degree of approximation are quadrics of the inaccuracies of the directly measurable variables $X_i$ ($i = 1, 2,..., n$). We prove that they describe a certain type of quadric hypersurfaces of parabolic class in the fields of the inaccuracies.

Our main results for the maximum inaccuracies in second degree of approximation can shortly be summarised as following:

1) A complete algebraic classification of the quadric hypersurfaces of the maximum inaccuracies of the indirectly measurable variable $Y$ is made;

2) A dimensionless scale for the quality (accuracy) evaluation of the experiment is introduced.

Besides the main results we researched in detail the two simplest partial cases in the algebraic classification and give example for the application of the dimensionless scale when experimental research commences.

**2. Preliminary results.** Let $Y = f(X_1, X_2,..., X_n)$ be an indirectly measurable variable depending on the directly measurable variables $X_i$ ($i = 1, 2,..., n$). Moreover, let us have $k_i$ number of observations of the directly measurable variable $X_i$ made in an experiment, which yield the values, respectively, $x_{i1}, x_{i2},..., x_{ik_i}$ ($i = 1, 2,..., n$).

The value of the partial derivative $\dfrac{\partial f}{\partial X_i}$, calculated on the $m$-th observation we denote with $\dfrac{\partial f}{\partial x_{im}}$ ($m = 1, 2,..., k_i$), and the arithmetic mean of the absolute values of the values of the partial derivative $\dfrac{\partial f}{\partial X_i}$ ($i = 1, 2,..., n$) we denote with $\overline{\left|\dfrac{\partial f}{\partial x_i}\right|} = \dfrac{1}{k_i}\sum_{m=1}^{k_i}\left|\dfrac{\partial f}{\partial x_{im}}\right|$. Analogously, we

denote with $\overline{\left|\dfrac{x_i}{f}\cdot\dfrac{\partial f}{\partial x_i}\right|} = \dfrac{1}{k_i}\sum_{m=1}^{k_i}\left|\dfrac{x_{im}}{f}\cdot\dfrac{\partial f}{\partial x_{im}}\right|$ the arithmetic mean of the absolute values of the values of $\dfrac{X_i}{f}\cdot\dfrac{\partial f}{\partial X_i}$ ($i=1,2,...,n$).

According to [5, 7] the maximum absolute inaccuracy of an indirectly measurable variable $Y$ in first degree of approximation is the linear function

(1) $$\Delta^1 Y = \Delta^1 f = \sum_{i=1}^{n}\overline{\left|\dfrac{\partial f}{\partial x_i}\right|}\cdot\left|\Delta X_i\right|$$

of the absolute values of the absolute inaccuracies $\Delta X_i$ of the directly measurable variables $X_i$ ($i=1,2,...,n$). According to [6, 7] the maximum relative inaccuracy of the indirectly measurable variable $Y$ in first degree of approximation is the linear function

(2) $$\dfrac{\Delta^1 Y}{Y} = \dfrac{\Delta^1 f}{|f|} = \sum_{i=1}^{n}\overline{\left|\dfrac{x_i}{f}\cdot\dfrac{\partial f}{\partial x_i}\right|}\cdot\left|\dfrac{\Delta X_i}{X_i}\right|$$

of the absolute values of the relative inaccuracies $\dfrac{\Delta X_i}{X_i}$ of the directly measurable variables $X_i$ ($i=1,2,...,n$).

For $i=1,2,...,n$ the coefficients $\overline{\left|\dfrac{\partial f}{\partial x_i}\right|}$ in (1) are called coefficients of influence of the absolute inaccuracies $\Delta X_i$ in $\Delta^1 Y$, and the coefficients $\overline{\left|\dfrac{X_i}{f}\cdot\dfrac{\partial f}{\partial X_i}\right|}$ in (2) – coefficients of influence of the relative inaccuracies $\dfrac{\Delta X_i}{X_i}$ in $\dfrac{\Delta^1 Y}{Y}$.

Our approach [5, 6, 7] is for $i=1,2,...,n$ to assume the coefficients of influence in (1) and (2) to be constants (within the given experiment), and the absolute inaccuracies $\Delta X_i$ in (1) and the relative inaccuracies $\dfrac{\Delta X_i}{X_i}$ in (2) to be variables.

If we look at $\Delta X_1, \Delta X_2,...,\Delta X_n, \Delta^1 Y$ as a system of generalised orthogonal coordinates [5, 7], then we get $n+1$-dimensional metric hyperspace $F^{n+1}$, in which (1) is an equation of a hyperplane passing through the beginning of the coordinate system. The hyperspace $F^{n+1}$ we call a space of the absolute inaccuracy of $Y$, and $\Delta^1 Y$ we call a plane of the absolute inaccuracy of $Y$. Analogically, $\dfrac{\Delta X_1}{X_1}, \dfrac{\Delta X_2}{X_2},...,\dfrac{\Delta X_n}{X_n}, \dfrac{\Delta^1 Y}{Y}$ is a system of generalised orthogonal coordinates [6, 7] of $n+1$-dimensional metric hyperspace $F_{n+1}$, in which (2) is an equation of a hyperplane passing through the beginning of the coordinate system. The hyperspace $F_{n+1}$ we call a space of the relative inaccuracy of $Y$, and $\frac{\Delta^1 Y}{Y}$ we call a plane of the relative inaccuracy of $Y$.

From (2) we get the equation $\alpha : \sum_{i=1}^{n} A_i \cdot \left| \frac{\Delta X_i}{X_i} \right| - \frac{\Delta^1 Y}{Y} = 0$ of the hyperplane $\alpha$ of the relative inaccuracy of $Y$, where $A_i = \left| \frac{x_i}{f} \cdot \frac{\partial f}{\partial x_i} \right| = const \geq 0$ ($i = 1, 2, ..., n$). Let us look into the hyperplane $\varepsilon : \frac{\Delta^1 Y}{Y} = 0$ also. Following [6, 7], we assume $\varepsilon$ to be a sample plane in the space of the relative inaccuracy in the sense that it corresponds to an imaginary ideal perfectly accurate experiment. Strictly speaking, such experiment is impossible and the plane $\varepsilon$ is unreachable. But by increasing the accuracy of the real experiment the plane $\alpha$ approaches $\varepsilon$.

Thus, the smaller the deviation of the plane $\alpha$ of the real experiment from the sample plane $\varepsilon$ of the ideal experiment, i.e. the smaller the angle between the two planes, the more accurate the experiment is. This angle is equal to the angle between the normal vectors $\vec{n}_\alpha (A_1, A_2, ..., A_n, -1)$ of the plane $\alpha$ and $\vec{n}_\varepsilon (0, 0, ..., 0, -1)$ of the plane $\varepsilon$. Then the value of the cosine

$$k_\alpha = \cos \angle (\vec{n}_\alpha, \vec{n}_\varepsilon) = \frac{1}{\sqrt{A_1^2 + A_2^2 + ... + A_n^2 + 1}}$$

of this angle is a coefficient of accuracy in a dimensionless scale [6, 7], i.e. a numerical characteristic of the quality of the experiment.

As $k_\alpha = \cos \angle (\vec{n}_\alpha, \vec{n}_\varepsilon)$, then the scale for evaluating the quality of the experiment is the interval $[0,1]$. An experiment is as accurate as the value of the coefficient of accuracy $k_\alpha$ is closer to 1, and is as inaccurate as the value of the coefficient of accuracy $k_\alpha$ is closer to 0. The value $k_\alpha = 1$ corresponds to the ideal perfectly accurate experiment and the value $k_\alpha = 0$ corresponds to the ideal absolutely inaccurate experiment.

**3. Complete algebraic classification of the maximum inaccuracies in second degree of approximation.** The variables defined correspondingly by

(3)
$$\Delta^2 Y = \Delta^2 f = \sum_{i,j=1}^{n} \overline{\left|\frac{\partial^2 f}{\partial x_i \partial x_j}\right|} \cdot |\Delta X_i| |\Delta X_j|$$

and

(4)
$$\frac{\Delta^2 Y}{Y} = \frac{\Delta^2 f}{|f|} = \sum_{i,j=1}^{n} \overline{\left|\frac{x_i x_j}{f} \cdot \frac{\partial^2 f}{\partial x_i \partial x_j}\right|} \cdot \left|\frac{\Delta X_i}{X_i}\right| \left|\frac{\Delta X_j}{X_j}\right|,$$

where the coefficients of influence $\overline{\left|\frac{\partial^2 f}{\partial x_i \partial x_j}\right|}$ and $\overline{\left|\frac{x_i x_j}{f} \cdot \frac{\partial^2 f}{\partial x_i \partial x_j}\right|}$ are constants, and $\Delta X_i$ and $\frac{\Delta X_i}{X_i}$ are variables ($i,j = 1,2,...,n$) we call, respectively, maximum absolute inaccuracy of second order and maximum relative inaccuracy of second order [7, 8].

The maximum absolute inaccuracy $\Delta Y$ of an indirectly measurable variable $Y$ in second degree of approximation, according to [8], is

(5)
$$\Delta Y = \Delta^1 Y + \frac{1}{2} \Delta^2 Y,$$

and the maximum relative inaccuracy $\frac{\Delta Y}{Y}$ of $Y$ in second degree of approximation is

(6)
$$\frac{\Delta Y}{Y} = \frac{\Delta^1 Y}{Y} + \frac{1}{2} \frac{\Delta^2 Y}{Y}.$$

Having in mind (1), (3) and (5) we get that the maximum absolute inaccuracy of the indirectly measurable variable $Y$ in second degree of approximation is presented in the following form

$$\Delta Y = \sum_{i=1}^{n} A_i \cdot |\Delta X_i| + \frac{1}{2} \sum_{i,j=1}^{n} A_{ij} \cdot |\Delta X_i| |\Delta X_j|,$$

where $A_i$ and $A_{ij}$ are non-negative constants, and $\Delta X_i$ are variables ($i,j = 1,2,...,n$). Analogically, from (2), (4) and (6) we get that the maximum relative inaccuracy of the indirectly measurable variable $Y$ in second degree of approximation is represented in the following form

$$\frac{\Delta Y}{Y} = \sum_{i=1}^{n} A_i \cdot \left|\frac{\Delta X_i}{X_i}\right| + \frac{1}{2} \sum_{i,j=1}^{n} A_{ij} \cdot \left|\frac{\Delta X_i}{X_i}\right| \left|\frac{\Delta X_j}{X_j}\right|,$$

where $A_i$ and $A_{ij}$ are non-negative constants, and $\frac{\Delta X_i}{X_i}$ are variables ($i,j = 1,2,...,n$).

Then in the corresponding spaces of the inaccuracies of $Y$ we can look into

$\Delta X_1, \Delta X_2, ..., \Delta X_n, \Delta Y$ and $\dfrac{\Delta X_1}{X_1}, \dfrac{\Delta X_2}{X_2}, ..., \dfrac{\Delta X_n}{X_n}, \dfrac{\Delta Y}{Y}$ as a system of generalised orthogonal coordinates, which shows that *the maximum inaccuracies of the indirectly measurable variable Y in second degree of approximation are quadrics of the inaccuracies of the directly measurable variables* $X_1, X_2, ..., X_n$. They describe quadric hypersurfaces passing through the beginning of the coordinate systems in these spaces.

We aim to produce a complete algebraic classification of the quadric hypersurfaces of the maximum inaccuracies of the indirectly measurable variable $Y$ in second degree of approximation. We will firstly look into the hypersurfaces of the maximum absolute inaccuracy and we will introduce some new notations for clarity.

With $y_1, y_2, ..., y_n, y_{n+1}$ we denote the system of generalised orthogonal coordinates $\Delta X_1, \Delta X_2, ..., \Delta X_n, \Delta Y$ in the space $F^{n+1}$. Then each such hypersurface can be represented in the form

$$y_{n+1} = \sum_{i=1}^{n} A_i y_i + \frac{1}{2} \sum_{i,j=1}^{n} A_{ij} y_i y_j,$$

where $A_i$ and $A_{ij}$ are non-negative constants ($i, j = 1, 2, ..., n$) and at least one of them is different from zero.

We construct a symmetric $(n+1) \times (n+1)$ matrix $A = (a_{ij})$ the following way: in the first $n$ rows and $n$ columns for $i, j = 1, 2, ..., n$ the elements $a_{ii} = \dfrac{1}{2} A_{ii}$ and $a_{ij} = a_{ji} = \dfrac{1}{4} A_{ij}$ are placed when $j \neq i$, and the elements in the $n+1$-th row and $n+1$-th column are zeros. Therefore, $A$ is the matrix of the quadratic part of the hypersurface. Moreover, let $b_i = \dfrac{1}{2} A_i$ for $i = 1, 2, ..., n$, $b_{n+1} = -1$, and $b^T = (b_1, b_2, ..., b_{n+1})$ is the transposed vector of $b$.

Therefore, the general form of the quadric hypersurfaces of the maximum absolute inaccuracies in second degree of approximation can be written in the form

(7) $$y^T A y + 2 b^T y = 0,$$

where $y^T = (y_1, y_2, ..., y_{n+1})$ is a vector (point) of the hypersurface [9, 10]. If we introduce the symmetric matrix $M = \begin{pmatrix} A & b \\ b^T & 0 \end{pmatrix}$ of type $(n+2) \times (n+2)$ and the $(n+2)$-dimensional vector $\hat{y} = \begin{pmatrix} y \\ 1 \end{pmatrix}$, then the equation of the hypersurface can also be written in the form $\hat{y}^T M \hat{y} = 0$.

We will also use the following basic facts [9, 10].

1) The characteristic roots $\lambda_1, \lambda_2, ..., \lambda_{n+1}$ of $A$ and the determinants $\det A$ and $\det M$ do not change their values when a linear transformation is applied (specifically – after translation and rotation with respect to the coordinate axis) of the hypersurface.

2) The quadratic hypersurface can be brought in a canonical form (unambiguously) with an orthogonal linear transformation. In other words, there is such an orthogonal matrix $U$, that the matrix $U^T A U$ is diagonal and on its diagonal there are the characteristic roots $\lambda_1, \lambda_2, ..., \lambda_{n+1}$ of $A$.

In order to produce a complete algebraic classification of the quadric hypersurfaces of the maximum absolute inaccuracy in second degree of approximation it is enough to produce such classification of their canonical types. Prior to describing the general case we will have a detailed look at the partial cases $n = 1$ and $n = 2$.

**Proposition 3.1.** *In the space $F^2$ the line of the maximum absolute inaccuracy of second degree of approximation can be only one of the following two kinds: a parabola and a straight line.*

**Proof.** In this case the quadric hypersurfaces of the maximum absolute inaccuracy in second degree of approximation are lines of second degree in the two dimensional plane $F^2$.

As in this case the calculations are simple we will find the canonical form of the line $a_{11} y_1^2 + b_1 y_1 - y_2 = 0$ directly. If $a_{11} = 0$, then the line is straight. If $a_{11} \neq 0$, then the linear transformation $y_1 = \frac{1}{\sqrt{a_{11}}} z_1$, $y_2 = \frac{b_1}{\sqrt{a_{11}}} z_1 + 2 p z_2$ brings the line to the parabola $z_1^2 = 2 p z_1$ [9].

**Proposition 3.2.** *In the space $F^3$ the surface of the maximum absolute inaccuracy in second degree of approximation can be only one of the following 4 kinds: an elliptic paraboloid, a hyperbolic paraboloid, a parabolic cylinder or a plane.*

**Proof.** In this case the quadric hypersurfaces of the maximum absolute inaccuracy in second degree of approximation are surfaces of second degree in the three dimensional space $F^3$. From (7) the surface of the maximum absolute inaccuracy in second degree of approximation has an equation $a_{11} y_1^2 + a_{22} y_2^2 + 2 a_{12} y_1 y_2 + 2 b_1 y_1 + 2 b_2 y_2 - y_3 = 0$, in which the coefficients are non-negative.

**3.2.1.** If the matrix $A = (a_{ij})$ is zero, then the surface is a plane. Depending on the number of non-zero coefficients $b_i$ ($i = 1, 2$) there are three kinds of such planes.

**3.2.2.** If $A$ is non-zero, then let $U$ be an orthogonal matrix which brings $A$ in a diagonal form with elements on its diagonal $\lambda_1, \lambda_2, \lambda_3$. As $A \neq 0$ and $\det A = 0$, then one or two of its characteristic roots are equal to 0.

**3.2.2.A.** Let $\lambda_1 \neq 0$, $\lambda_2 \neq 0$ and $\lambda_3 = 0$. We use the rotation $y = U y'$. Then the surface is represented by the equation $\lambda_1(y_1')^2 + \lambda_2(y_2')^2 + 2c_1 y_1' + 2c_2 y_2' + 2c_3 y_3' = 0$, where $(c_1, c_2, c_3) = b^T U$. The presence of terms of second degree with respect to $y_1'$ and $y_2'$ allows us to remove the terms of first degree with respect to $y_1'$ and $y_2'$ by the substitution $y_1'' = y_1' + \frac{c_1}{\lambda_1}$, $y_2'' = y_2' + \frac{c_2}{\lambda_2}$. Thus, we get the equation $\lambda_1(y_1'')^2 + \lambda_2(y_2'')^2 + 2c_3 y_3' + q = 0$, where $q = -\frac{c_1^2}{\lambda_1} - \frac{c_2^2}{\lambda_2}$. Let the matrix $M$ transform to the matrix $N$ after these transformations. As $\det M = -\frac{1}{4}\lambda_1\lambda_2$, and from the last equation we have $\det N = -c_3^2 \lambda_1 \lambda_2$, then $c_3^2 = \frac{1}{4}$, i.e. $c_3 \neq 0$. Further, there are two possible cases for the characteristic roots of the matrix $A$.

**3.2.2.A.a.** Let $\lambda_1\lambda_2 > 0$. With the translation $y_3'' = y_3' + \frac{q}{2c_3}$ with respect to the axis $y_3'$ we remove the constant term in the equation and get $\lambda_1(y_1'')^2 + \lambda_2(y_2'')^2 + 2c_3 y_3'' = 0$. We divide the two sides of the equation with $-c_3$ and get the canonical form $\frac{(y_1'')^2}{p_1} + \frac{(y_2'')^2}{p_2} = 2y_3''$ of elliptic paraboloid [9].

**3.2.2.A.b.** Let $\lambda_1\lambda_2 < 0$. With analogical transformations we get the canonical form $\frac{(y_1'')^2}{p_1} - \frac{(y_2'')^2}{p_2} = 2y_3''$ of hyperbolic paraboloid [9].

**3.2.2.B.** Let $\lambda_2 = \lambda_3 = 0$ and $\lambda_1 \neq 0$. After the rotation $y = Uy'$ we get the equation $\lambda_1(y_1')^2 + 2c_1 y_1' + 2c_2 y_2' + 2c_3 y_3' = 0$, where $(c_1, c_2, c_3) = b^T U$ and $c = \sqrt{c_2^2 + c_3^2} \neq 0$. With translation with respect to the axis $y_1'$ we remove the term of first degree with respect to $y_1'$ and get $\lambda_1(y_1'')^2 + 2c_2 y_2' + 2c_3 y_3' + q = 0$. After the rotation $y_2'' = \frac{c_2}{c} y_2' + \frac{c_3}{c} y_3'$, $y_3'' = \frac{c_2}{c} y_2' - \frac{c_3}{c} y_3'$ we get $\lambda_1(y_1'')^2 + 2c y_2'' + q = 0$. With translation with respect to the axis $y_2''$ we remove the constant term and by dividing by $\lambda_1$ we reach the canonical form $(y_1'')^2 = 2p_2 y_2''$ of parabolic cylinder [9].

**3.3. Note.** We would like to remind about some facts regarding the general classification of the hypersurfaces of second degree in an Euclidean space with at least 2 dimensions. Let $A$ be the matrix of the quadratic part of the hypersurface. Depending on the values of the characteristic roots of the matrix $A$ these hypersurfaces are divided into three classes – elliptic, hyperbolic and parabolic.

The elliptic class is characterised by the condition that all the characteristic roots of $A$ are non-zero and are all positive or all negative. This means that $A$ is either positive or negative definite (the Sylvester criterion gives a necessary and sufficient condition for this) [10, 11].

The hyperbolic class is characterised by the condition that all the characteristic roots of $A$ are non-zero but at least two of them have different signs. This means that $A$ is not defined with a certain sign [10, 11].

The parabolic class is characterised by the condition that at least one of the characteristic roots of $A$ is zero, i.e. $\det A = 0$.

**Theorem 3.4.** *In the space $F^{n+1}$, when $n > 2$, the surface of the maximum absolute inaccuracy in second degree of approximation can be only one of the following 4 kinds: an elliptic paraboloid, a hyperbolic paraboloid, a parabolic cylinder or a hyperplane.*

**Proof.** Let us look at a quadric hypersurface of a maximum absolute inaccuracy in second degree of approximation. Formula (7) yields that the hypersurface has an equation $\sum_{i,j=1}^{n} a_{ij} y_i y_j + \sum_{i=1}^{n+1} b_i y_i = 0$, in which the coefficients are non-negative.

**3.4.1.** If $A = 0$, then the hypersurface is a hyperplane. Depending on the number of non-zero coefficients $b_i$ there are $n+1$ kinds of such hyperplanes.

**3.4.2.** If $A \neq 0$, then as $(n+1)$-th row and $(n+1)$-th column of $A$ are zero, then $\det A = 0$, i.e. this hypersurface can only be of parabolic class.

From $A \neq 0$ and $\det A = 0$ it follows that it is not possible all the characteristic roots of $A$ to be equal to zero. If the rank of $A$ is equal to $r$, where $1 \leq r \leq n$, then the number of non-zero characteristic roots of $A$ is also $r$. Then, if $U$ is an orthogonal matrix which brings $A$ in a diagonal form, then $U^T A U$ is a diagonal matrix on the diagonal of which the characteristic roots $\lambda_1, \lambda_2, \ldots, \lambda_r, 0, \ldots, 0$ of $A$ are.

After the rotation $y = Uy'$ the equation of the hypersurface assumes the form $\lambda_1 (y'_1)^2 + \lambda_2 (y'_2)^2 + \ldots + \lambda_r (y'_r)^2 + 2c_1 y'_1 + \ldots + 2c_{n+1} y'_{n+1} = 0$, where $(c_1, c_2, \ldots, c_{n+1}) = b^T U$. With the translation $y''_i = y'_i + \dfrac{c_i}{\lambda_i}$, $1 \leq i \leq r$, we remove the linear terms containing $y'_1, y'_1, \ldots, y'_r$. Then the equation becomes $\lambda_1 (y''_1)^2 + \lambda_2 (y''_2)^2 + \ldots + \lambda_r (y''_r)^2 + 2c_{r+1} y'_{r+1} + \ldots + 2c_{n+1} y'_{n+1} + q = 0$, where $q = -\dfrac{c_1^2}{\lambda_1} - \ldots - \dfrac{c_r^2}{\lambda_r}$. The type of the hypersurface depends on the value of $r$ and thus we will address the following two cases.

**3.4.2.A.** Let $r = n$. Then the equation has the form $\lambda_1 (y''_1)^2 + \lambda_2 (y''_2)^2 + \ldots + \lambda_r (y''_n)^2 + 2c_{n+1} y'_{n+1} + q = 0$. Let after the translations above are applied the matrix $M$ be transformed into the matrix $N$. As $\det M = -\dfrac{1}{4} \lambda_1 \lambda_2 \ldots \lambda_n$, and from the last equation we have $\det N = -c_{n+1}^2 \lambda_1 \lambda_2 ... \lambda_n$, then $c_{n+1}^2 = \dfrac{1}{4}$, i.e. $c_{n+1} \neq 0$. Further on there are two possible subcases.

**3.4.2.A.a.** Let all the non-zero characteristic roots $\lambda_1, \lambda_2, ..., \lambda_n$ have the same sign. With the translation $y''_{n+1} = y'_{n+1} + \dfrac{q}{2c_{n+1}}$ we remove the constant term in the equation. We divide the two sides of the equation by $-c_{n+1}$ and get $\dfrac{z_1^2}{p_1} + ... + \dfrac{z_n^2}{p_n} = 2z_{n+1}$, where all the parameters $p_1, p_2, ..., p_n$ are with the same sign. Thus, the canonical form of the quadric hypersurface is an elliptic paraboloid.

**3.4.2.A.b.** Let at least two of the characteristic roots $\lambda_1, \lambda_2, ..., \lambda_n$ of $A$ be with different signs. With analogical transformations we get the equation $\dfrac{z_1^2}{p_1} + ... + \dfrac{z_k^2}{p_k} - \dfrac{z_{k+1}^2}{p_{k+1}} - ... - \dfrac{z_n^2}{p_n} = 2z_{n+1}$, where all the parameters $p_1, p_2, ..., p_n$ are with the same sign. Therefore, the canonical form of the quadric hypersurface is a hyperbolic paraboloid. Depending on the count of the positive signs there are $n-1$ kinds of such paraboloids. We will say that the hyperbolic paraboloid is of type $k$, $1 \leq k \leq n-1$, if it has exactly $k$ positive coefficients.

**3.4.2.B.** Let $r < n$, then the researched surface has the equation $\lambda_1 (y''_1)^2 + \lambda_2 (y''_2)^2 + ... + \lambda_r (y''_r)^2 + 2c_{r+1} y'_{r+1} + ... + 2c_{n+1} y'_{n+1} + q = 0$, and $c = \sqrt{c_{r+1}^2 + ... + c_{n+1}^2} \neq 0$. In $n+1-r$-dimensional subspace $y''_1 = y''_2 = ... = y''_r = 0$ we will apply such a rotation that the hypersurface has no projections on the axes $y'_{r+2}, ..., y'_{n+1}$. In order to do that we construct an orthogonal matrix $V$ of type $(n+1-r) \times (n+1-r)$ with first row $\left(\dfrac{c_{r+1}}{c}, ..., \dfrac{c_{n+1}}{c}\right)$. Then after the rotation $(y''_{r+1}, ..., y''_{n+1})^T = V(y'_{r+1}, ..., y'_{n+1})^T$ we get the equation $\lambda_1 (y''_1)^2 + \lambda_2 (y''_2)^2 + ... + \lambda_r (y''_r)^2 + 2c_{r+1} y''_{r+1} + q'' = 0$. With translation with respect to the axis $y''_{r+1}$ we remove the constant term, we divide by $-c_{r+1}$ and get the equation $\dfrac{z_1^2}{p_1} + ... + \dfrac{z_r^2}{p_r} = 2z_{r+1}$, where all parameters $p_1, p_2, ..., p_n$ are with the same sign. This is the canonical form of the quadric hypersurface and it is a parabolic cylinder. Depending on the rank of $A$ there are $n-1$ kinds of such cylinders. We will say that the parabolic cylinder is of type $r$, $1 \leq r \leq n-1$, if the rank of $A$ is equal to $r$.

**3.5. Note.** The research of quadric hypersurfaces of the maximum relative inaccuracy follows exactly the same scheme and reaches the same conclusions as for the absolute. It is enough only with $y_1, y_2, ..., y_n, y_{n+1}$ to denote the system of generalised orthogonal coordinates $\dfrac{\Delta X_1}{X_1}, \dfrac{\Delta X_2}{X_2}, ..., \dfrac{\Delta X_n}{X_n}, \dfrac{\Delta Y}{Y}$ in the space $F_{n+1}$.

**4. Dimensionless scale for evaluating the quality (accuracy) of the experiment when determining the maximum inaccuracies in second degree of approximation.** Following the notation from Section 3 and Note 3.5 the equation of quadric hypersurface of the maximum relative inaccuracy of the indirectly measurable variable $Y$ in second degree of approximation can be written in the form $\alpha : \sum_{i,j=1}^{n} a_{ij} y_i y_j + \sum_{i=1}^{n} b_i y_i - y_{n+1} = 0$. Here $a_{ij}$ and $b_i$ are constants ($i, j = 1, 2, ..., n$) and without loss of generality we can assume that at least one of the coefficients $a_{ij}$ is non-zero. Let us look at the quadric hypersurface $\varepsilon : y_{n+1} = 0$ with coefficients $a_{ij} = b_i = 0$. Analogically to the idea of sample plane, given in Section 2, we can assume $\varepsilon$ to be sample hypersurface in the space of the relative inaccuracy, corresponding to an imaginary idea perfectly accurate experiment. Once again, the deviation of the hypersurface $\alpha$ of the real experiment from the sample hypersurface $\varepsilon$ of the ideal experiment can be used as a measurement for the accuracy of the experiment.

In the case of hyperplanes, described in Section 2, the angle between the normal vectors $\vec{n_\alpha}$ of $\alpha$ and $\vec{n_\varepsilon}$ of $\varepsilon$ is the same regardless of the choice of the point in which $\vec{n_\alpha}$ is applied. The problem with the hypersurfaces is that in each point of $\alpha$ the angle between $\vec{n_\alpha}$, applied in this point, and $\vec{n_\varepsilon}$ is different. This problem can be overcame as following.

Let us denote with $\nabla \alpha = \left( \dfrac{\partial \alpha}{\partial y_1}, \dfrac{\partial \alpha}{\partial y_2}, ..., \dfrac{\partial \alpha}{\partial y_n}, \dfrac{\partial \alpha}{\partial y_{n+1}} \right)$ the gradient of $\alpha$. For each point $P = (p_1, p_2, ..., p_{n+1})$ of $\alpha$ the vector $\nabla \alpha (P)$ is perpendicular to $\alpha$ in this point [12]. Let us denote with $\overline{p_i}$ the arithmetic mean of all the values of $\dfrac{\partial \alpha}{\partial y_i}$ given an experiment and $\overline{P} = (\overline{p_1}, \overline{p_2}, ..., \overline{p_n})$. Then we assume that the deviation of the hypersurface $\alpha$ of the real experiment from the sample hypersurface $\varepsilon$ of the ideal experiment is the angle between the vectors $\overline{\nabla \alpha} = \nabla \alpha(\overline{P})$ and $\vec{n_\varepsilon}$. More precisely, the value of the cosine

(8) $$k_\alpha = \cos \angle (\overline{\nabla \alpha}, \vec{n_\varepsilon}) = \dfrac{1}{\sqrt{(\overline{p_1})^2 + (\overline{p_2})^2 + ... + (\overline{p_n})^2 + 1}}$$

we assume to be a coefficient of accuracy in a dimensionless scale, i.e. for numerical characteristics of the quality of the experiment. Thus, again, the scale for determining the quality of the experiment is the interval $[0,1]$. Moreover, an experiment is as accurate as the value of the coefficient of accuracy $k_\alpha$ is closer to 1, and is as accurate as the value of the coefficient of accuracy $k_\alpha$ is closer to 0. (The value $k_\alpha = 1$ corresponds to the ideal perfectly accurate experiment, and the value $k_\alpha = 0$ corresponds to the ideal absolute inaccurate experiment.)

From the formula for $k_\alpha$ we can give the following criteria for accuracy of an experiment:

*An experiment is most accurate if and only if the sum of the squares of the coefficients $\left(\overline{p_1}\right)^2 + \left(\overline{p_2}\right)^2 + ... + \left(\overline{p_n}\right)^2$ assumes its lowest value.*

**5. Example for determining the canonical form of a hypersurface of the maximum inaccuracy in second degree of approximation and evaluating the quality of an experiment.** The viscosity $\eta$ of a tested liquid with known density $\rho$ (given fixed temperature) can be determined using an Ostwald viscometer. In order to do that a control liquid with known viscosity $\eta_0$ and density $\rho_0$ (given the same temperature) is chosen. The times for outflowing $t_0$ and $t$ from the viscometer of the same volume of control liquid and the test liquid are measured. The viscosity $\eta$ is calculated by the formula $\eta = \eta_0 \dfrac{\rho t}{\rho_0 t_0}$.

We determined by this method the viscosity of 40% water solution of glycerine given a temperature of $20°C$ with control liquid distilled water which given this temperature has a viscosity $\eta_0 = 1.002 \times 10^{-3}$ **Pa.s** and density $\rho_0 = 998.23$ **kg/m³** [13]. Initially, we made 5 measurements of the time $t_0$ of outflowing of the distilled water with volume 10 ml. Afterwards we measured 5 times the time $t$ of outflowing of 40% water solution of glycerine with the same volume, which given $20°C$ has density $\rho = 1098.4$ **kg/m³** [13]. The results are shown in Table 1.

| № | 1 | 2 | 3 | 4 | 5 |
|---|---|---|---|---|---|
| $t_0$ [s] | 40.2 | 38.5 | 38.9 | 39.2 | 39.7 |
| $t$ [s] | 11.6 | 11.7 | 11.8 | 12 | 12.1 |

Table 1. Time of outflowing of distilled water and 40% water solution of glycerine from an Ostwald viscometer

According to (2) the maximum relative inaccuracy of the indirectly measurable variable $\eta$ in first degree of approximation is $\dfrac{\Delta^1 \eta}{\eta} = \dfrac{\Delta t}{t} + \dfrac{\Delta t_0}{t_0}$, and according to (4) the maximum relative inaccuracy of second order is $\dfrac{\Delta^2 \eta}{\eta} = \dfrac{\Delta t}{t} \cdot \dfrac{\Delta t_0}{t_o} + 2\left(\dfrac{\Delta t_0}{t_0}\right)^2$. Then following (6) the maximum relative inaccuracy $\dfrac{\Delta \eta}{\eta}$ of $\eta$ in second degree of approximation is $\dfrac{\Delta \eta}{\eta} = \dfrac{\Delta t}{t} + \dfrac{\Delta t_0}{t_0} + \dfrac{1}{2} \dfrac{\Delta t}{t} \cdot \dfrac{\Delta t_0}{t_0} + \left(\dfrac{\Delta t_0}{t_0}\right)^2$. As one of the characteristic roots of the matrix of the quadratic form is zero and the other two are with different signs ($\lambda_1 = \frac{2-\sqrt{5}}{4}, \lambda_2 = \frac{2+\sqrt{5}}{4}$), then the canonical form of the surface of the maximum relative inaccuracy $\frac{\Delta \eta}{\eta}$ in second degree of approximation is a hyperbolic paraboloid.

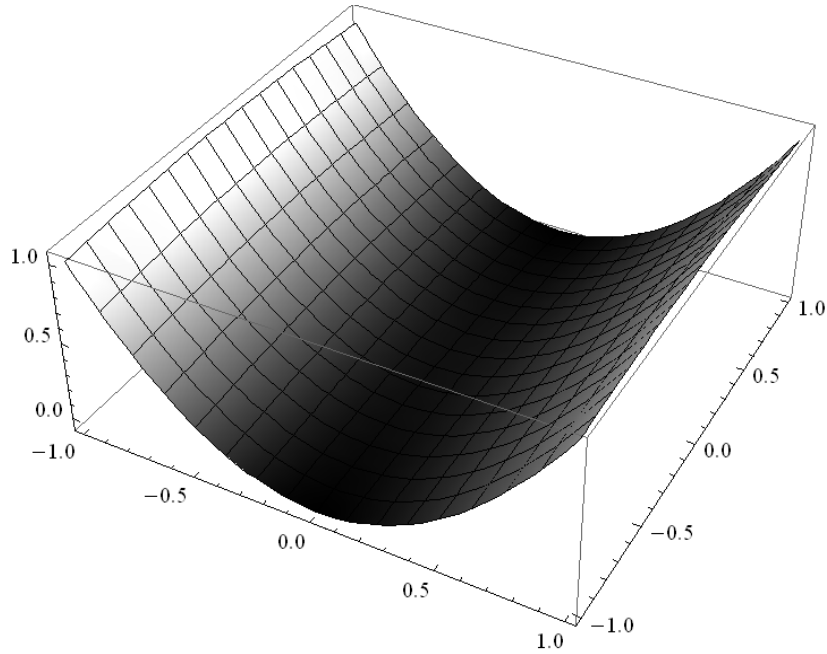

Figure 1. The canonical form of the hyperbolic paraboloid $\frac{\Delta \eta}{\eta}$

In order to calculate the coefficient of accuracy when determining the maximum relative inaccuracy $\frac{\Delta \eta}{\eta}$ in second degree of approximation we firstly find the gradient $\nabla \eta = \left(1 + \frac{1}{2} \cdot \frac{\overline{\Delta t_0}}{t_0}, 1 + \frac{1}{2} \cdot \frac{\overline{\Delta t}}{t} + 2\frac{\overline{\Delta t_0}}{t_0}, -1\right)$. By calculating its value given $\frac{\overline{\Delta t_0}}{t_0} = 0.013$ and $\frac{\overline{\Delta t}}{t} = 0.014$, we find $\overline{\nabla \eta} = (1.006, 1.033, -1)$. Then, according to (8), the coefficient of accuracy is $k_\eta = 0.57$ and as it is close to 1 the experiment is sufficiently accurate.

**6. Discussion.** The parabolic class hypersurfaces is characterised with the property that it divides the hypersurfaces from the elliptic and hyperbolic class. Namely, let us look at a random hypersurface of second degree of $n+1$ variables in the $n+1$-dimensional real space $R^{n+1}$ and $A$ is the symmetric matrix of its quadratic part. Let $m = \frac{(n+1)(n+2)}{2}$ and let us equate $A$ with the vector $a = (a_{11}, a_{21}, ..., a_{n+1,1}, a_{22}, a_{32}, ..., a_{n+1,2}, ..., a_{n+1,n+1})$ from $m$-dimensional real space $R^m$. Then if $A$ describes the parabolic class hypersurfaces in $R^{n+1}$ then its correspondent vector $a$ describes a cone in $R^m$ with an apex in the beginning of the coordinate system. This cone divides $R^m$ to two parts: internal (corresponding to the elliptic class in $R^{n+1}$) and external (corresponding to the hyperbolic class in $R^{n+1}$).

**7. Conclusion.** An analytical approach for determining the maximum inaccuracies of an indirectly measurable variable is developed from the interval analysis by constructing an interval extension for the given function [14]. The idea of the differential interval extension, which is a type of interval extension, is conceptually close to the proposed in the paper algebraic approach.

Our method for representation of the maximum inaccuracies of an indirectly measurable variable is universal, because it is applicable to different scientific fields, different mathematical models, different kinds of measurements, different types of input data.

Finally, we would like to point out that the suggested by us algebraic classification of the maximum inaccuracies of an indirectly measurable variable $Y$ in second degree of approximation is complete, thus of a particular importance. According to it these maximum inaccuracies are quadrics of the inaccuracies of $X_1, X_2, ..., X_n$ and describe exactly certain kinds quadric hypersurfaces of parabolic class.

**Aknowledgements**

The paper is partially funded by the "Scientific Studies" Fund of the Plovdiv University 'P. Hilendarski' as per the contract NI15-FMIIT-004/24.04.2015.